\documentclass[a4paper,11pt]{article}

\input{diagrams}

\newcommand{\cl}[1]{\ensuremath{\mathcal {#1}}}

\newcommand{\lra}{\ensuremath{\longrightarrow}}
\newcommand{\map}[1]{\ensuremath{\stackrel{{#1}}{\lra}}}

\newcommand{\nspan}[5]
{\setlength{\unitlength}{1mm}
\begin{center}
\begin{picture}(30,30)(15,15)      %
\put(17,23){\makebox(0,0)[t]{${#1}$}}  
\put(32,37){\makebox(0,0)[b]{${#3}$}}  
\put(45,23){\makebox(0,0)[t]{${#5}$}}  

\put(29,35){\vector(-1,-1){10}}  %
\put(33,35){\vector(1,-1){10}}  %
\put(23,32){\makebox(0,0)[r]{${#2}$}} 
\put(39,32){\makebox(0,0)[l]{${#4}$}}

\end{picture}
\end{center}}

\newcommand{\bpoint}[1]{
\begin{itemize}
\item {#1}	
\end{itemize}}

\newtheorem{theorem}{Theorem}[section]

\newtheorem{definition}[theorem]{Definition}

\newcommand{\numroman}{\renewcommand{\labelenumi}{\roman{enumi})}}

\newcommand{\pica}{\begin{center} \input}
\newcommand{\picz}{\end{center}}

\newlength{\leng}
\newlength{\fontleng}
\newcommand{\sunit}{\setlength{\unitlength}{1mm}}
\newcommand{\setleng}[1]{\setlength{\leng}{#1}
    \setlength{\unitlength}{\leng}}
\newcommand{\setunit}[1]{\setlength{\unitlength}{#1}}

\newcommand{\diaga}[3]{

\put(#1,#2){\setlength{\leng}{#3}
\setlength{\unitlength}{1.3\leng} %

\put(1,2){\line(1,0){5}} %
\put(1,3){\line(1,0){5}} %
\put(1,4){\line(1,0){5}} %
\put(1,5){\line(1,0){5}} %
\put(1,6){\line(1,0){5}} %
\put(4,0){\line(1,1){4}} %
\put(4,8){\line(1,-1){4}} }}

\newcommand{\diagb}[3]{

\put(#1,#2){\setlength{\leng}{#3}
\setlength{\unitlength}{1.3\leng} %

\put(3,1){\line(0,1){5}} %
\put(4,1){\line(0,1){5}} %
\put(5,1){\line(0,1){5}} %
\put(6,1){\line(0,1){5}} %
\put(7,1){\line(0,1){5}} %
\put(1,4){\line(1,1){4}} %
\put(5,8){\line(1,-1){4}}}}

\newcommand{\diagc}[3]{

\put(#1,#2){
\setlength{\unitlength}{#3} %

\put(1,2){\line(1,0){5}} %
\put(1,3){\line(1,0){5}} %
\put(1,4){\line(1,0){5}} %
\put(1,1){\line(1,0){5}} %
\put(3.5,-1.5){\line(1,1){4}} %
\put(3.5,6.5){\line(1,-1){4}}}}

\newcommand{\diagd}[3]{

\put(#1,#2){\setlength{\leng}{#3}
\setlength{\unitlength}{0.7\leng} %

\put(0,6){\line(1,-2){3}} %
\put(2,7){\line(1,-2){3}} %
\put(4,8){\line(1,-2){3}} %
\put(6,9){\line(1,-2){3}} %
\put(0,0){\line(1,0){7}} %
\put(7,0){\line(2,3){4}}}}

\newcommand{\diage}[3]{

\put(#1,#2){\setlength{\leng}{#3}
\setlength{\unitlength}{0.7\leng} %

\put(0,-6){\line(1,2){3}} %
\put(2,-7){\line(1,2){3}} %
\put(4,-8){\line(1,2){3}} %
\put(6,-9){\line(1,2){3}} %
\put(0,0){\line(1,0){7}} %
\put(7,0){\line(2,-3){4}} }}

\newcommand{\diagf}[3]{

\put(#1,#2){
\setlength{\unitlength}{#3} %

\put(0,0){\line(1,0){4}}  %
\put(0,0){\line(2,3){2}}  %
\put(4,0){\line(-2,3){2}}}}

\newcommand{\diagg}[3]{

\put(#1,#2){
\setlength{\unitlength}{#3} %

\put(0,0){\line(1,0){4}}  %
\put(0,0){\line(2,3){2}}  %
\put(4,0){\line(-2,3){2}} %
\qbezier(0,0)(2,1)(4,0) }}

\newcommand{\diagi}[3]{

\put(#1,#2){
\setlength{\unitlength}{#3} %

\put(0,0){\line(1,0){4}}  %
\put(0,0){\line(2,3){2}}  %
\put(4,0){\line(-2,3){2}} %
\qbezier(0,0)(2,1)(4,0)   %
\qbezier(0,0)(2,2)(4,0) %
\qbezier(0,0)(2,3)(4,0) %
}}

\newcommand{\diagj}[3]{

\put(#1,#2){
\setlength{\unitlength}{#3} %

\put(0,0){\line(1,0){4}}  %
\qbezier(0,0)(2,2)(4,0) }}

\newcommand{\diagk}[3]{

\put(#1,#2){
\setlength{\unitlength}{#3} %

\put(0,0){\line(1,0){4}}  %
\qbezier(0,0)(2,1)(4,0)   %
\qbezier(0,0)(2,2)(4,0) %
}}

\newcommand{\diagl}[3]{

\put(#1,#2){
\setlength{\unitlength}{#3} %

\put(0,0){\line(1,0){4}}  %
\qbezier(0,0)(2,1)(4,0)   %
\qbezier(0,0)(2,2)(4,0) %
\qbezier(0,0)(2,3)(4,0) %
}}

\newcommand{\diagm}[3]{

\put(#1,#2){
\setlength{\unitlength}{#3} %

\put(0,0){\line(1,0){4}}  %
\put(0,0){\line(2,3){2}}  %
\put(4,0){\line(-2,3){2}} %
\qbezier(2,3)(4,2)(4,0)}}

\newcommand{\diago}[3]{

\put(#1,#2){
\setlength{\unitlength}{#3} %

\put(0,0){\line(1,0){4}}  %
\put(0,0){\line(2,3){2}}  %
\put(4,0){\line(-2,3){2}} %
\qbezier(2,3)(4,2)(4,0)  %
\qbezier(2,3)(4.9,2.6)(4,0)  %
\qbezier(2,3)(5.6,3.4)(4,0)}}

\newcommand{\diagp}[3]{

\put(#1,#2){

\setlength{\unitlength}{#3}  %
\fontsize{#3}{20}

\put(0,0){\makebox(0,0){$\equiv$}}  %
\put(0.35,0){\makebox(0,0){$\rangle$}}}}

\newcommand{\diagq}[3]{

\put(#1,#2){

\setlength{\unitlength}{#3}  %

\qbezier(0.5,0)(2.5,0.5)(3,4)   %
\put(3,4){\makebox(0,0)[bc]{{\tiny $\Delta$}}} }}

\newcommand{\diagr}[3]{

\put(#1,#2){\setlength{\leng}{#3}

\setlength{\unitlength}{1.4\leng}

\qbezier(0.5,0)(2.5,1)(3,4)   %
\put(0.5,-0){\makebox(0,0)[c]{{\tiny $\Delta$}}} }}

\newcommand{\diags}[3]{

\put(#1,#2){\setlength{\leng}{#3}
\setlength{\unitlength}{2\leng} %

\diagj{0}{0}{2\leng}     %
\diagj{7}{0}{2\leng}        %
\diagp{5.5}{0.5}{3\leng}}}

\newcommand{\diagt}[3]{

\put(#1,#2){\setlength{\leng}{#3}
\setlength{\unitlength}{2\leng} %

\diagk{0}{0}{2\leng}       %
\diagj{7}{0}{2\leng}         %
\diagp{5.5}{0.5}{3\leng} }}

\newcommand{\diagu}[3]{

\put(#1,#2){\setlength{\leng}{#3}
\setlength{\unitlength}{2\leng} %

\diagl{0}{0}{2\leng}       %
\diagj{7}{0}{2\leng}         %
\diagp{5.5}{0.5}{3\leng} }}

\newcommand{\diagv}[3]{

\put(#1,#2){\setlength{\leng}{#3}
\setlength{\unitlength}{2\leng} %

\diagf{0}{0}{2\leng}       %
\diagf{7}{0}{2\leng}         %
\diagp{5.5}{1.5}{3\leng} }}

\newcommand{\diagw}[3]{

\put(#1,#2){\setlength{\leng}{#3}
\setlength{\unitlength}{2\leng} %

\diagm{0}{0}{2\leng}       %
\diagf{7}{0}{2\leng}         %
\diagp{5.5}{1.5}{3\leng} }}

\newcommand{\diagy}[3]{

\put(#1,#2){\setlength{\leng}{#3}
\setlength{\unitlength}{2\leng} %

\diago{0}{0}{2\leng}       %
\diagf{7.4}{0}{2\leng}         %
\diagp{6.4}{1.5}{3\leng} }}

\newcommand{\diagz}[3]{

\put(#1,#2){\setlength{\leng}{#3}
\setlength{\unitlength}{2\leng} %

\diagg{0}{0}{2\leng}       %
\diagf{7}{0}{2\leng}         %
\diagp{5.5}{1.5}{3\leng} }}

\newcommand{\diagab}[3]{

\put(#1,#2){\setlength{\leng}{#3}
\setlength{\unitlength}{2\leng} %

\diagi{0}{0}{2\leng}       %
\diagf{7}{0}{2\leng}         %
\diagp{5.5}{1.5}{3\leng} }}

\newcommand{\diagad}[3]{

\put(#1,#2){

\setlength{\unitlength}{#3}  %

\qbezier(0.5,0)(2.5,-0.5)(3,-4)   %
\put(3,-4){\makebox(0,0)[t]{{\tiny $\nabla$}}}}}

\newcommand{\diagac}[3]{

\put(#1,#2){\setlength{\leng}{#3}
\setlength{\unitlength}{2\leng} %

\diags{0}{0}{\leng}       %
\diagt{0}{6}{\leng}         %
\diagt{9}{3}{\leng}    %
\diagq{9}{0.5}{\leng} %
\diagad{9}{6.2}{\leng}}}

\newcommand{\diagae}[3]{

\put(#1,#2){\setlength{\leng}{#3}
\setlength{\unitlength}{2\leng} %

\diagac{0}{0}{\leng}       %
\diagt{0}{6}{\leng}         %
\diagu{28}{3}{\leng}    %
\diagc{22}{2.5}{\leng}}} %

\newcommand{\diagaf}[3]{

\put(#1,#2){\setlength{\leng}{#3}
\setlength{\unitlength}{2\leng} %

\diagv{0}{0}{\leng}       %
\diagv{17}{0}{\leng}         %
\diagc{12}{0}{\leng}}} %

\newcommand{\diagag}[3]{

\put(#1,#2){\setlength{\leng}{#3}
\setlength{\unitlength}{2\leng} %

\diagw{0}{0}{\leng}       %
\diagy{17}{0}{\leng}         %
\diagc{12}{0}{\leng}    %
\diagu{-4.5}{5}{\leng}    %
\diagr{3}{2.5}{\leng}}} %

\newcommand{\diagah}[3]{

\put(#1,#2){\setlength{\leng}{#3}
\setlength{\unitlength}{2\leng} %

\diagv{0}{0}{\leng}       %
\diagy{17}{0}{\leng}         %
\diagc{12}{0}{\leng}    %
\diagy{-5}{4.5}{\leng}    %
\diagr{2}{2.2}{\leng}}} %

\newcommand{\diagai}[3]{

\put(#1,#2){\setlength{\leng}{#3}
\setlength{\unitlength}{2\leng} %

\diagv{0}{0}{\leng}       %
\diagab{17}{0}{\leng}         %
\diagc{12}{0}{\leng}    %
\diagab{-5}{4.5}{\leng}    %
\diagr{2}{2.2}{\leng}}} %

\newcommand{\diagaj}[3]{

\put(#1,#2){\setlength{\leng}{#3}
\setlength{\unitlength}{2\leng} %

\diagz{0}{0}{\leng}       %
\diagab{17}{0}{\leng}         %
\diagc{12}{0}{\leng}    %
\diagu{-7.5}{-3}{\leng}    %
\diagq{1}{-2.5}{\leng}}} %

\newcommand{\diagak}[3]{

\put(#1,#2){
\setlength{\unitlength}{#3} %

\put(0,0){\line(1,0){4}}  %
\put(0,0){\line(2,3){2}}  %
\put(4,0){\line(-2,3){2}} %
\qbezier(2,3)(4.9,2.6)(4,0)}}

\newcommand{\diagal}[3]{

\put(#1,#2){
\setlength{\unitlength}{#3} %

\put(0,0){\line(1,0){4}}  %
\put(0,0){\line(2,3){2}}  %
\put(4,0){\line(-2,3){2}} %
\qbezier(0,0)(2,2)(4,0) }}

\newcommand{\diagam}[3]{

\put(#1,#2){\setlength{\leng}{#3}
\setlength{\unitlength}{2\leng} %

\diagak{0}{0}{2\leng}       %
\diagf{7.3}{0}{2\leng}         %
\diagp{6}{1.5}{3\leng} }}

\newcommand{\diagan}[3]{

\put(#1,#2){\setlength{\leng}{#3}
\setlength{\unitlength}{2\leng} %

\diagal{0}{0}{2\leng}       %
\diagf{7.3}{0}{2\leng}         %
\diagp{6}{1.5}{3\leng} }}

\newcommand{\diagas}[3]{

\put(#1,#2){\setlength{\leng}{#3}
\setlength{\unitlength}{\leng} %

\diagag{-15}{0}{\leng}       %
\diagae{-35}{25}{\leng}       %
\diaga{45}{10}{\leng}       %
\diagw{75}{7}{\leng}       %
\diagc{100}{7}{\leng}       %
\diagy{110}{7}{\leng}       %
\diagac{57}{20}{\leng}       %
\diagy{85}{-30}{\leng}       %
\diagd{75}{-12}{\leng}       %
\diage{115}{-7}{\leng}       %
\diagb{90}{-15}{\leng}       %

\qbezier(94,27)(92,15)(83.3,11.7)  %
\put(84,11.5){\makebox(0,0)[r]{{\tiny $\Delta$}}}

\qbezier(32.5,32)(29,26)(12,24)  %
\qbezier(12,24)(-9,21)(-13,13)  %
\put(-13,13){\makebox(0,0)[t]{{\scriptsize $\nabla$}}}}}

\newcommand{\diagat}[3]{

\put(#1,#2){\setlength{\leng}{#3}
\setlength{\unitlength}{\leng} %

\diagah{-15}{25}{\leng}       %
\diagaf{-30}{0}{\leng}       %
\diaga{45}{10}{\leng}       %
\diagah{75}{7}{\leng}       %
\diagy{85}{-30}{\leng}       %
\diagd{75}{-12}{\leng}       %
\diage{115}{-7}{\leng}       %
\diagb{90}{-15}{\leng}       %

\qbezier(14.5,4)(14,8)(5,13)  %
\qbezier(5,13)(-4,19)(-4.5,26)  %
\put(-4,26){\makebox(0,0)[b]{{\scriptsize $\Delta$}}}}}

\newcommand{\diagau}[3]{

\put(#1,#2){\setlength{\leng}{#3}
\setlength{\unitlength}{\leng} %

\diagai{-15}{25}{\leng}       %
\diagaf{-30}{0}{\leng}       %
\diaga{45}{10}{\leng}       %
\diagai{75}{7}{\leng}       %
\diagab{85}{-30}{\leng}       %
\diagd{75}{-12}{\leng}       %
\diage{115}{-7}{\leng}       %
\diagb{90}{-15}{\leng}       %

\qbezier(14.5,4)(14,8)(5,13)  %
\qbezier(5,13)(-4,19)(-4.5,26)  %
\put(-4,26){\makebox(0,0)[b]{{\scriptsize $\Delta$}}}}}

\newcommand{\diagav}[3]{

\put(#1,#2){\setlength{\leng}{#3}
\setlength{\unitlength}{\leng} %

\diagaj{-15}{30}{\leng}       %
\diagae{-35}{-5}{\leng}       %
\diaga{45}{10}{\leng}       %
\diagz{75}{25}{\leng}       %
\diagc{100}{25}{\leng}       %
\diagab{110}{25}{\leng}       %
\diagac{57}{5}{\leng}       %
\diagab{85}{-30}{\leng}       %
\diagd{75}{-12}{\leng}       %
\diage{115}{-7}{\leng}       %
\diagb{90}{-15}{\leng}      %

\qbezier(79,24)(80,20)(86,18)  %
\qbezier(86,18)(92,16)(93,12)   %
\put(79,24){\makebox(0,0)[b]{{\tiny $\Delta$}}}

\qbezier(-19,23)(-18,15)(5,12)  %
\qbezier(5,12)(30,10)(31.5,3)  %
\put(-19,23){\makebox(0,0)[b]{{\scriptsize $\Delta$}}}}}

\newcommand{\diagaw}[3]{
\put(#1,#2){\setlength{\leng}{#3}
\setlength{\unitlength}{4\leng} %
\fontsize{3.8\leng}{15}

\put(0,0){\line(1,0){4}}  %
\put(0,0){\line(2,3){2}}  %
\put(4,0){\line(-2,3){2}}
\put(0.3,1.6){$f$}   %
\put(3.4,1.6){$g$}   %
\put(2,-1){$\bar{g}$} %
\put(1.8,1){$\Downarrow$}

}}

\newcommand{\diagax}[3]{

\put(#1,#2){\setlength{\leng}{#3}
\setlength{\unitlength}{4\leng} %
\fontsize{3.8\leng}{15}

\diagam{0}{0}{2\leng}  %
\put(0.4,1.6){$f$}          %
\put(3.3,1.5){$\phi$}          %
\put(5.5,2.8){$u_{31}$} }}

\newcommand{\diagay}[3]{

\put(#1,#2){\setlength{\leng}{#3}
\setlength{\unitlength}{4\leng} %
\fontsize{3.8\leng}{15}

\diagan{0}{0}{2\leng}  %
\put(0.4,1.6){$f$}          %
\put(1.7,0.5){{\small $\bar\phi$}}          %
\put(5.5,2.8){$u_{32}$} }}

\newcommand{\diagaz}[3]{

\put(#1,#2){\setlength{\leng}{#3}
\setlength{\unitlength}{\leng} %
\fontsize{3.8\leng}{15}

\diagay{0}{0}{\leng}  %
\diagax{0}{20}{\leng} %
\put(36.5,15){$\|$}     %
\put(36,4){$\phi_1$}   %
\put(36,24){$\phi_1$}  %
\put(7,6){$u_2$}  %
\put(7,24){$u_2'$} }}

\newcommand{\diagba}[3]{

\put(#1,#2){\setlength{\leng}{#3}
\setlength{\unitlength}{\leng} %
\fontsize{3.8\leng}{15}

\diagag{0}{0}{\leng}  %
\put(1,14){$\phi$}     %
\put(25,8){$u_{41}$}   %
\put(45,7){$\phi_1$}   %
\put(10,-2){$u_{31}$}     %

}}

\newcommand{\diagbb}[3]{

\put(#1,#2){\setlength{\leng}{#3}
\setlength{\unitlength}{\leng} %
\fontsize{3.8\leng}{15}

\diagah{0}{0}{\leng}  %
\put(10,-2){$\phi_2$}     %
\put(25,8){$u_{42}$}   %
\put(45,7){$\phi_1$}   %

   }}

\newcommand{\diagbc}[3]{

\put(#1,#2){\setlength{\leng}{#3}
\setlength{\unitlength}{\leng} %
\fontsize{3.8\leng}{15}

\diagai{0}{0}{\leng}  %
\put(10,-2){$\phi_2$}     %
\put(25,8){$u_{43}$}   %
\put(43,7){$\phi_3$}   %
}}

\newcommand{\diagbd}[3]{

\put(#1,#2){\setlength{\leng}{#3}
\setlength{\unitlength}{\leng} %
\fontsize{3.8\leng}{15}

\diagaj{0}{0}{\leng}  %
\put(-5,-11){$\bar\phi$}     %
\put(25,8){$u_{44}$}   %
\put(43,7){$\phi_3$}   %
\put(10,7){$u_{32}$}     %
}}

\newcommand{\diagbf}[3]{

\put(#1,#2){\setlength{\leng}{#3}
\setlength{\unitlength}{\leng} %
\fontsize{3.8\leng}{15}

\diagas{0}{0}{\leng}  %
\put(12,40){$\phi$} %
\put(10,-6){$u_{41}$} %
\put(48,23){$u_{51}$} %
\put(103,17){$\phi_1$} %
\put(123,-15){$\phi_2$} %
\put(103,-12){$u_{52}$}  %

}}

\newcommand{\diagbg}[3]{

\put(#1,#2){\setlength{\leng}{#3}
\setlength{\unitlength}{\leng} %
\fontsize{3.8\leng}{15}

\diagat{0}{0}{\leng}  %
\put(-4,-7){$\phi_4$} %
\put(11,33){$u_{42}$} %
\put(48,23){$u_{54}$} %
\put(101,17){$\phi_3$} %
\put(123,-15){$\phi_2$} %
\put(103,-12){$u_{53}$}  %

}}

\newcommand{\diagbh}[3]{

\put(#1,#2){\setlength{\leng}{#3}
\setlength{\unitlength}{\leng} %
\fontsize{3.8\leng}{15}

\diagau{0}{0}{\leng}  %
\put(-4,-7){$\phi_4$} %
\put(11,33){$u_{43}$} %
\put(48,23){$u_{55}$} %
\put(101,17){$\phi_5$} %
\put(123,-15){$\phi_6$} %
\put(103,-12){$u_{56}$}  %
}}

\newcommand{\diagbi}[3]{

\put(#1,#2){\setlength{\leng}{#3}
\setlength{\unitlength}{\leng} %
\fontsize{3.8\leng}{15}

\diagav{0}{0}{\leng}  %
\put(10,40){$u_{44}$} %
\put(11,-6){$\bar\phi$} %
\put(48,23){$u_{58}$} %
\put(103,35){$\phi_7$} %
\put(123,-15){$\phi_6$} %
\put(103,-12){$u_{57}$}  %
}}

\sunit

\newcommand{\onetwo}[8]{

\put(#1,#2){ \setlength{\leng}{#3} %
\setlength{\unitlength}{\leng} %

\diagj{0}{0}{6\leng} %
\put(-1,0){\makebox(0,0)[r]{#4}} %
\put(26,0){\makebox(0,0)[l]{#5}} %
\put(12,7){\makebox(0,0)[b]{#6}} %
\put(12,-2){\makebox(0,0)[t]{#7}} %
\put(12,2.5){\makebox(0,0)[c]{$\Downarrow$}}  %
\put(14,2.5){\makebox(0,0)[l]{#8}} }}

\sunit

\sunit

\sunit

\sunit

\newcommand{\three}[4]{
\put(#1,#2){ \setlength{\leng}{#3} %
\setlength{\unitlength}{\leng} %
\diagp{0}{0}{6\leng}  %
\put(0,4){\makebox(0,0)[b]{#4}} }}

\sunit

\sunit

\sunit

\newcommand{\twotwo}[7]{

\put(0,0){ \setlength{\leng}{1mm} %
\setlength{\unitlength}{\leng} %
\diagf{0}{0}{4\leng}  %
\put(-1,-1){\makebox(0,0)[tr]{#1}}   
\put(8,14){\makebox(0,0)[b]{#2}}  %
\put(18,-1){\makebox(0,0)[tl]{#3}}  %

\put(3,6){\makebox(0,0)[br]{#4}}   
\put(14,6){\makebox(0,0)[bl]{#5}}  %
\put(8,-2){\makebox(0,0)[t]{#6}}  %

\put(8,5){\makebox(0,0)[c]{$\Downarrow$}}  %
\put(9,5){\makebox(0,0)[l]{#7}}  %
}}

\sunit

\newcommand{\twotwob}[7]{

\put(0,0){ \setlength{\leng}{1mm} %
\setlength{\unitlength}{\leng} %
\diagf{0}{0}{4\leng}  %
\put(-1,-1){\makebox(0,0)[tr]{#1}}   
\put(8,14){\makebox(0,0)[b]{#2}}  %
\put(18,-1){\makebox(0,0)[tl]{#3}}  %

\put(3,6){\makebox(0,0)[br]{#4}}   
\put(14,6){\makebox(0,0)[bl]{#5}}  %
\put(8,-2){\makebox(0,0)[t]{#6}}  %

\put(8,5){\makebox(0,0)[c]{#7}}  %

}}

\newcommand{\diagbr}[9]{

\put(0,0){\setlength{\leng}{1mm}
\setlength{\unitlength}{\leng} %

\diagal{0}{0}{4\leng}
\put(-1,-1){\makebox(0,0)[tr]{#1}}   
\put(8,14){\makebox(0,0)[b]{#2}}  %
\put(18,-1){\makebox(0,0)[tl]{#3}}  %

\put(3,6){\makebox(0,0)[br]{#4}}   
\put(14,6){\makebox(0,0)[bl]{#5}}  %
\put(8,-2){\makebox(0,0)[t]{#7}}  %

\put(8,7){\makebox(0,0)[c]{#8}}  %

\put(10,3.3){\makebox(0,0)[bl]{#6}}  %
\put(8,1.6){\makebox(0,0)[c]{#9}}  %
}}

\newcommand{\diagbs}[3]{

\put(#1,#2){
\setlength{\unitlength}{#3} %

\put(0,0){\line(1,0){4}}  %
\put(0,0){\line(2,3){2}}  %
\put(4,0){\line(-2,3){2}} %
\qbezier(2,3)(4.9,2.6)(4,0)}}

\newcommand{\diagbt}[9]{

\put(0,0){\setlength{\leng}{1mm}
\setlength{\unitlength}{\leng} %

\diagbs{0}{0}{4\leng}
\put(-1,-1){\makebox(0,0)[tr]{#1}}   
\put(8,14){\makebox(0,0)[b]{#2}}  %
\put(18,-1){\makebox(0,0)[tl]{#3}}  %

\put(3,6){\makebox(0,0)[br]{#4}}   
\put(12,4){\makebox(0,0)[c]{#6}}  %
\put(8,-2){\makebox(0,0)[t]{#7}}  %

\put(8,6){\makebox(0,0)[c]{#9}}  %

\put(15,6){\makebox(0,0)[c]{#8}}  %
\put(17,9){\makebox(0,0)[bl]{#5}}  %
}}

\newcommand{\scalecq}{

\put(0,0){
\setlength{\unitlength}{0.3\leng} %

\put(20,90){\line(0,1){20}}      %
\put(40,90){\line(0,1){20}}      %
\put(80,90){\line(0,1){20}}      %

\put(20,90){\line(1,0){60}}      %
\put(20,90){\line(1,-1){30}}     %
\put(80,90){\line(-1,-1){30}}    %
\put(50,40){\line(0,1){20}}      %
\put(60,104){\makebox[0pt]{$\cdots$}} }}

\newcommand{\scalecr}[5]{

\put(0,0){
\setunit{0.3\leng} %
\scalecq
\put(20,114){\makebox[0pt]{#1}}   %
\put(40,114){\makebox[0pt]{#2}}   %
\put(80,114){\makebox[0pt]{#3}}   %
\put(50,76){\makebox(0,0){#5}}      
\put(45,38){\makebox(0,0)[t]{#4}}  }} 

\newcommand{\scalecs}[9]{

\put(0,0){
\setunit{0.3\leng} %
\put(20,0){ \scalecq
\put(20,114){\makebox[0pt]{#4}}   %
\put(40,114){\makebox[0pt][r]{#5}}   %
\put(80,114){\makebox[0pt]{#6}}   %
\put(50,76){\makebox(0,0){#9}}      
\put(50,36){\makebox(0,0)[t]{#7}}} \setunit{0.2\leng}
\put(39.8,124){\setleng{0.67\leng} \scalecr{#1}{#2}{#3}{}{#8}}}}

\newcommand{\scalect}[4]{
\setlength{\unitlength}{0.3\leng} %

\put(20,0){\line(0,1){20}} %
\put(20,20){\line(-2,3){20}} %
\put(20,20){\line(2,3){20}}  %
\put(0,50){\line(1,0){40}}  %
\put(0,50){\line(0,1){20}}  %
\put(40,50){\line(0,1){20}} %

\put(15,0){\makebox(0,0)[r]{#3}} %
\put(0,75){\makebox(0,0)[b]{#1}}  %
\put(40,75){\makebox(0,0)[b]{#2}}  %
\put(20,40){\makebox(0,0)[c]{#4}}}

\newcommand{\scalecu}[3]{
\setlength{\unitlength}{0.3\leng} %

\put(20,0){\line(0,1){20}} %
\put(20,20){\line(-2,3){20}} %
\put(20,20){\line(2,3){20}}  %
\put(0,50){\line(1,0){40}}  %
\put(20,50){\line(0,1){20}} %

\put(15,0){\makebox(0,0)[r]{#2}} %
\put(20,75){\makebox(0,0)[b]{#1}}  %
\put(20,40){\makebox(0,0)[c]{#3}}}

\newcommand{\scalecv}[6]{
\put(3,0){\scalect{}{#3}{#4}{#6}} %
\put(-1.5,21){\setleng{0.75\leng} \scalecu{#1}{#2}{#5}}}

\newcommand{\scalecy}[7]{
\put(-9,0){\scalecr{#1}{#2}{#3}{#4}{#6}} %
\put(1.5,-2.5){\setleng{0.75\leng} \scalecu{}{#5}{#7}}}

\usepackage{amsfonts}
\usepackage{amssymb}
\usepackage{amsmath}
\usepackage{fancyheadings}
\usepackage{latexcad}

\begin{document}

\title{An alternative characterisation of universal cells in opetopic
$n$-categories}
\author{Eugenia Cheng\\ \\Department of Pure Mathematics, University
of Cambridge\\E-mail: e.cheng@dpmms.cam.ac.uk}
\date{October 2002}
\maketitle

\begin{abstract}

We address the fact that composition in an opetopic weak $n$-category is in
general not unique and hence is not a well-defined operation.  We define
composition with a given $k$-cell in an $n$-category by a span of
$(n-k)$-categories.  We characterise such a cell as universal if its composition
span gives an equivalence of $(n-k)$-categories.  

\end{abstract}

\setcounter{tocdepth}{3}
\tableofcontents

\numroman

\section*{Introduction}
\addcontentsline{toc}{section}{Introduction}

In this paper we give a characterisation of universality in the opetopic theory
of $n$-categories.  

The opetopic definition of $n$-category proceeds in two stages.  First, the
language for describing $k$-cells is constructed.  This is the theory of
opetopes.  Then, a concept of universality is introduced, to deal with
composition and coherence.  Eventually, we have the following definition:

\begin{definition} An {\em opetopic (weak) $n$-category} is an opetopic set in
which 

\begin{enumerate}

\item Every niche has an $n$-universal occupant.
\item Every composite of $n$-universals is $n$-universal.

\end{enumerate}
\end{definition}

The word `composite' is used in the following sense: given any universal cell,
its target cell is said to be a composite of its source cells.  Thus in the
opetopic theory, composites are not necessarily unique.  

In \cite{che7} and \cite{che8} we examine the various approaches to the theory
of opetopes (\cite{bd1}, \cite{hmp1}, \cite{lei1}) and prove that they are
equivalent.  In the present paper we turn our attention to the second stage of
the definition of opetopic $n$-category, concerning universality.  

There are many ways of characterising universal cells, just as there are many
ways of characterising, say, isomorphisms in a category.  The original
definition given by Baez and Dolan generalises the result `$f$ is an
isomorphism if and only if any morphism with the same domain factors through it
uniquely'. The characterisation proposed here is motivated by another familiar
result in categories, that $f$ is an isomorphism if and only if composition with
$f$ is an isomorphism.  

In a category, ``composition with $f$'' is a function on homsets; however, composition in
an opetopic $n$-category is not uniquely defined, that is, $\_ \circ f$ is not
a well-defined operation.  One way of dealing with this would be to {\em
choose} composites in order to make $\_ \circ f$ into a well-defined
operation.  However, this is not in the spirit of the opetopic definition.  To
avoid making such choices we instead definte ``composition with $f$'' as a span
of hom-$(n-k)$-categories.  This ``composition span'' gives all possible ways of
composing with $f$.  

We begin in Section~\ref{defuniv} by recalling the definition of universal
cells in an opetopic set $X$, given in \cite{che10}.  This is our motivation
for our new characteristion of universal cells in an $n$-category.  In
Section~\ref{char} we show how the composition span is used in this
characterisation, and in Section~\ref{univcon} we give the actual construction
of the composition span.  While we intend that the two notions of universality
should coincide when $X$ is an $n$-category, we do not include a proof here as
we currently lack an effective method for calculating in arbitrary dimensions. 
As a gesture towards this result, we include some low-dimensional examples in
Section~\ref{univdraw}, and in Section~\ref{comp} we prove that the notions do
indeed coincide for $n \leq 2$.  Finally in Section~\ref{conc} we make some
brief concluding remarks.

\label{altuniv}

\subsubsection*{Terminology and Notation}

\begin{enumerate}

\item Since we are concerned chiefly with weak $n$-categories we continue our
previous practice of omitting the word ``weak'' in general.

\item In this paper we will avoid any detailed discussion of the language of
multicategories and construction of opetopic sets; this has been discussed in
detail in our earlier work (\cite{che7}, \cite{che8}, \cite{che9}).  We will
adopt the (more practical) method of Hermida, Makkai and Power (\cite{hmp1}),
picking one ordering of source elements in order to represent a symmetry class.
A $k$-cell has as its source a pasting diagram of $(k-1)$-cells, and as its
target a single $(k-1)$-cell.  For a general $k$-cell we write its source as
$\underline{a}$, say, to indicate a formal composite whose constituent
$(k-1)$-cells may be placed in some order.

\item Furthermore, we may adopt the following convention for 2-ary
cells.  A 2-ary $k$-cell $\alpha$ has the form

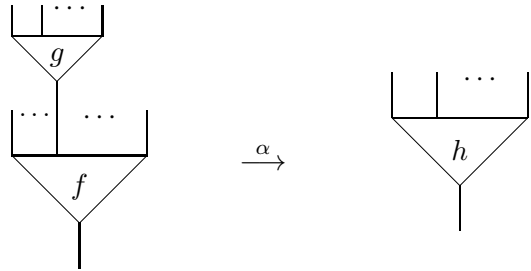
\begin{center} \setleng{1mm}
\begin{picture}(18,30)(8,8)
\scalecs{}{}{}{}{}{}{}{$g$}{$f$} \put(13,32.2){{\small $\ldots$}}
\end{picture}
\begin{picture}(42,18) \put(15,17){$\map{\alpha}$} \end{picture}
\begin{picture}(18,45)(14,3)
\scalecr{}{}{}{}{$h$}
\end{picture}
\end{center}

\noindent where $f$, $g$, and $h$ are $(k-1)$-cells (and
necessarily $k \geq 2$).  We write this $k$-cell as \[\alpha :
(f,g) \lra h\] employing this ordering of the source elements to
indicate that $f$ and $g$ are pasted at the target of $g$; we also
write $s_1(\alpha) = f$ and $s_2(\alpha) = g$.

\end{enumerate}

\bigskip {\bfseries Acknowledgements}

This work was supported by a PhD grant from EPSRC.  I would like to thank
Martin Hyland and Tom Leinster for their support and guidance.

\section{Definition of universality}
\label{defuniv}

We begin by recalling from \cite{che10} the definition of a universal cell in an
opetopic set.

Let $X$ be an opetopic set.  

\begin{definition} A $k$-cell $\alpha$ is {\em $n$-universal} if either
$k>n$ and $\alpha$ is unique in its niche, or $k \leq n$ and (1)
and (2) below are satisfied:

\renewcommand{\labelenumi}{(\arabic{enumi})}
\begin{enumerate}
\item Given any $k$-cell $\gamma$ in the same niche as $\alpha$, there is a
factorisation $u:(\beta,\alpha) \lra \gamma$

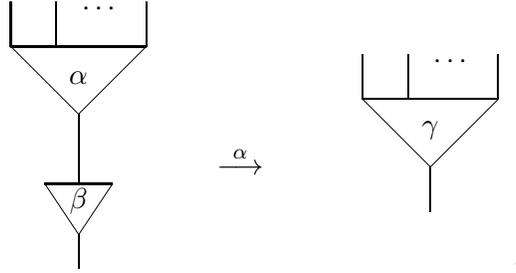
\begin{center}
\setleng{1mm}
\begin{picture}(18,40)(0,-2)
\scalecy{}{}{}{}{}{$\alpha$}{$\beta$}
\end{picture}
\begin{picture}(25,18) \put(5,12){$\map{\alpha}$} \end{picture}
\begin{picture}(18,30)(8,5)
\scalecr{}{}{}{}{$\gamma$}
\end{picture}.
\end{center}

\item Any such factorisation is $n$-universal.

\end{enumerate}
\end{definition}

\begin{definition} A factorisation $u:(b,a) \lra c$ of $k$-cells
is {\em $n$-universal} if $k>n$, or $k\leq n$ and (1) and (2)
below are satisfied:

\renewcommand{\labelenumi}{(\arabic{enumi})}
\begin{enumerate}
\item Given any $k$-cell $b'$ in the same frame as $b$, and any
$(k+1)$-cell \[v:(b',a) \lra c\] with $b'$ and $a$ pasted in the
same configuration as $b$ and $a$ in the source of $u$, there is a
factorisation of $(k+1)$-cells $(u,y)\lra v$

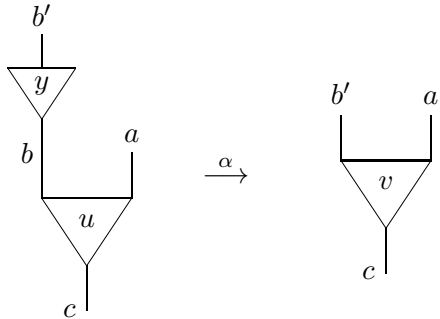
\begin{center}
\setleng{1mm}
\begin{picture}(18,45)(0,0)
\scalecv{$b'$}{$b$}{$a$}{$c$}{$y$}{$u$}
\end{picture}
\begin{picture}(30,18) \put(5,17){$\map{\alpha}$} \end{picture}
\begin{picture}(18,30)(8,-5)
\scalect{$b'$}{$a$}{$c$}{$v$}
\end{picture}
\end{center}

\item Any such factorisation is itself $n$-universal.

\end{enumerate}
\end{definition}

\section{An alternative characterisation}
\label{char}

We now examine the motivating example in categories.  Let
\cl{C} be a category and $f:A \longrightarrow B$ a morphism in
\cl{C}. Then we have a natural transformation
    \[H^f: \cl{C}(B,\_\ ) \lra \cl{C}(A,\_\ )\]
with components
    \[\_ \circ f : \cl{C}(B,C) \lra \cl{C} (A,C)\]
for each $C \in \cl{C}$.  Then
\[\begin{array}{rcl}
    f \mbox{ is an isomorphism} & \iff & H^f \mbox{ is an isomorphism} \\
    & \iff & \forall\  C \in \cl{C},\  \_ \circ f  \mbox{ is an isomorphism} \\
    & \iff & \mbox{ ``composition with $f$ is an isomorphism''}
    \end{array}\]

\noindent Here ``composition with $f$'' is a function on homsets.

Now let $X$ be an opetopic $n$-category and $f:\underline{a} \lra
b$ a $k$-cell in $X$.  Then given any $(k-1)$-cell $c$ we have
$(n-k)$-categories $X(b,c)$ and $X(\underline{a},c)$ whose 0-cells
are $k$-cells of $X$ with the appropriate source and target, and
whose $j$-cells are $(k+j)$-cells.

Since composition in an opetopic $n$-category is not uniquely
defined, we cannot expect $\_\circ f$ to be a well-defined
operation $X(b,c) \lra X(\underline{a},c)$. Instead, we will have
a span of $(n-k)$-categories

\nspan{X(b,c)}{\sigma_f}{C_f}{\tau_f}{X(\underline{a},c)}

\noindent where $C_f$ gives all possible ways of composing with
$f$. Here $\sigma_f$ and $\tau_f$ are $(n-k)$-functors i.e.\
morphisms of the underlying opetopic sets.  ($\sigma$ has more
properties that we will not discuss here.)

We then have the following characterisation.

\begin{definition}\label{univdef}  A $k$-cell $f$ is {\em universal} iff
    \begin{enumerate}
    \item $k>n$ and $f$ is unique in its niche, or
    \item $k \leq n$ and $\tau_f$ is an $(n-k)$-equivalence of
    $(n-k)$-categories.
    \end{enumerate}
\end{definition}

\begin{definition} An $m$-functor is an {\em $m$-equivalence of
$m$-categories} iff
    \begin{enumerate}
    \item it is an $(m-1)$-equivalence on hom-$(m-1)$-categories
    \item it is ``essentially surjective on 0-cells'' i.e.\
    surjective up to universal 1-cells
    \end{enumerate}
\end{definition}

We observe (without giving details) that since $\sigma_f$ will be
an $(n-k)$-equivalence of $(n-k)$-categories, the above condition
for universality will also result in $X(b,c)$ and
$X(\underline{a},c)$ being $(n-k)$-equivalent.

Furthermore it will follow from the construction of the
composition span that in an $n$-category the above definition is
equivalent to demanding ``on the nose'' surjectivity.  i.e.\ $f$
is universal iff $\forall \underline{x},y \in C_f$
    \[\tau:C_f(\underline{x},y) \lra X(\underline{\tau x},\tau
    y)\]
is surjective on objects.   This is a consequence of the fact that
composites of universals are universal in an opetopic
$n$-category.

In the next section we construct the composition span itself.

\section{Construction of composition span}
\label{univcon}

In this section we give the construction of a composition span; in
the next section we give some explicit examples at low dimensions.

Composition of $k$-cells is given by universal $(k+1)$-cells, so
in order to construct a composition span for a $k$-cell $f$, we
must assume that  for all $m>k$ the universal $m$-cells have been
defined.

We seek to construct a span of opetopic sets

\nspan{X(b,c)}{\sigma_f}{C_f}{\tau_f}{X(\underline{a},c).}

\noindent For convenience we write $C_f = C$, $\sigma_f = \sigma$
and $\tau_f = \tau$ . Also put $X(b,c) = X_1$ and
$X(\underline{a},c) = X_2$. Recall that a morphism $F:A \lra B$ of
opetopic sets has for each $j\geq 0$ a function
    \[F_j: A(j) \lra B(j)\]
such that, for each $j\geq 1$ a certain square

\begin{center}
\setlength{\unitlength}{0.2em}
\begin{picture}(40,48)(8,-1)      %

\put(10,10){\makebox(0,0){$.$}}  
\put(10,35){\makebox(0,0){$A(j)$}}  
\put(45,35){\makebox(0,0){$B(j)$}}  
\put(45,10){\makebox(0,0){$.$}}  

\put(17,35){\vector(1,0){21}}  
\put(17,10){\vector(1,0){21}}  
\put(10,30){\vector(0,-1){15}} 
\put(45,30){\vector(0,-1){15}} 

\put(8,23){\makebox(0,0)[r]{}} 
\put(47,23){\makebox(0,0)[l]{}} 
\put(27,37){\makebox(0,0)[b]{$F_j$}} 
\put(27,8){\makebox(0,0)[t]{$$}} 


\end{picture}
\end{center}

\noindent commutes, ensuring that ``underlying shapes are
preserved''.  So we seek for each $j\geq 0$ functions

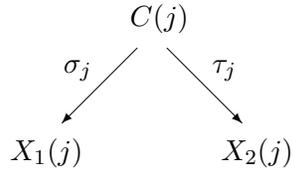
\nspan{X_1(j)}{\sigma_j}{C(j)}{\tau_j}{X_2(j)}

\noindent such that for each $j\geq 1$ a certain diagram

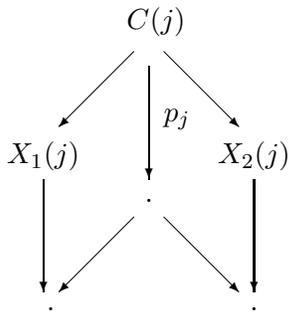
\begin{center} \setlength{\unitlength}{1mm}
\begin{picture}(35,55)(15,15)
\put(0,22){

\begin{picture}(80,30)      %
\put(17,23){\makebox(0,0)[t]{$X_1(j)$}}  
\put(32,37){\makebox(0,0)[b]{$C(j)$}}  
\put(45,23){\makebox(0,0)[t]{$X_2(j)$}}  

\put(29,35){\vector(-1,-1){10}}  %
\put(33,35){\vector(1,-1){10}}  %
\put(23,32){\makebox(0,0)[t]{$$}} 
\put(39,32){\makebox(0,0)[t]{$$}}

\put(17,18){\vector(0,-1){15}}  %
\put(31,33){\vector(0,-1){15}}  %
\put(45,18){\vector(0,-1){15}}  %
\put(33,26){\makebox(0,0)[l]{$p_j$}}

\end{picture}}

\put(0,0){
\begin{picture}(80,30)     %
\put(18,23){\makebox(0,0)[t]{$.$}}  
\put(31,37){\makebox(0,0)[b]{$.$}}  
\put(45,23){\makebox(0,0)[t]{$.$}}  

\put(29,35){\vector(-1,-1){10}}  %
\put(33,35){\vector(1,-1){10}}  %
\put(23,32){\makebox(0,0)[t]{$$}} 
\put(39,32){\makebox(0,0)[t]{$$}}
\end{picture}}

\end{picture}
\end{center}

\noindent commutes.  Then a $j$-cell $\theta \in C(j)$ exhibits
$\tau_j(\theta)\in X_2(j)$ as a composite of $f$ with
$\sigma_j(\theta) \in X_1(j)$. $p_j$ gives the frame of each
$j$-cell in $C$.

\begin{itemize}
\item $j=0$
\end{itemize}

Put \[C(0) = \{u \in \cl{U}(k+1)\  |\ s_2(u) = f \}\] where
$\cl{U}(m)$ is the set of 2-ary universal $m$-cells.  Put
$\sigma_0=s_1$ and $\tau_0=t$.

\begin{itemize}
\item$j=1$
\end{itemize}

A 1-frame in $C$ has the form $u_1 \lra u_2$.  We form the set of
occupants of this frame as follows.  Write
    \[\begin{array}{c} \cl{U}_1 = \{u \in \cl{U}(k+2)\ |\ s_1(u)=u_2\}\\
    \cl{U}_2 = \{u \in \cl{U}(k+2)\ |\ s_2(u)=u_1\} \end{array}\]
and form the pullback

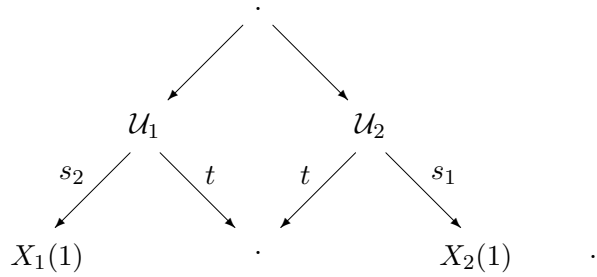
\begin{center}
\setlength{\unitlength}{1mm}
\begin{picture}(65,50)(15,10)
\put(15,17){

\begin{picture}(80,30)      %
\put(17,23){\makebox(0,0)[t]{$$}}  
\put(31,37){\makebox(0,0)[b]{$.$}}  
\put(45,23){\makebox(0,0)[t]{$$}}  

\put(29,35){\vector(-1,-1){10}}  %
\put(33,35){\vector(1,-1){10}}  %
\put(23,32){\makebox(0,0)[t]{$$}} 
\put(39,32){\makebox(0,0)[t]{$$}}

\end{picture}}

\put(0,0){
\begin{picture}(80,30)     %
\put(18,23){\makebox(0,0)[t]{$X_1(1)$}}  
\put(31,37){\makebox(0,0)[b]{$\cl{U}_1$}}  
\put(46,22){\makebox(0,0)[t]{$.$}}  

\put(29,35){\vector(-1,-1){10}}  %
\put(33,35){\vector(1,-1){10}}  %
\put(23,32){\makebox(0,0)[r]{$s_2$}} 
\put(39,32){\makebox(0,0)[l]{$t$}}
\end{picture}}

\put(30,0){
\begin{picture}(80,30)     %
\put(18,23){\makebox(0,0)[t]{$$}}  
\put(31,37){\makebox(0,0)[b]{$\cl{U}_2$}}  
\put(45,23){\makebox(0,0)[t]{$X_2(1)$}}  

\put(29,35){\vector(-1,-1){10}}  %
\put(33,35){\vector(1,-1){10}}  %
\put(23,32){\makebox(0,0)[r]{$t$}} 
\put(39,32){\makebox(0,0)[l]{$s_1$}}
\put(60,21){.}
\end{picture}}

\end{picture}
\end{center}

\begin{itemize}
\item $j>1$
\end{itemize}

For higher values of $j$ we construct for each $j$ a pullback over
$2^j$ subsets of $\cl{U}(k+j+1)$ as follows.

Let $\theta$ be a $(j-1)$-frame in $C$ with target $\alpha$.
$\alpha$ is a $(j-1)$-cell of $C$ so is a string of $2^{j-1}$
universal $(k+j)$-cells $u_1, \ldots, u_{j-1}$, say.  Now write
$\cl{U}=\cl{U}(k+j+1)$ and for each $1 \leq i \leq 2^{j-1}$
    \[\cl{U}_i=\{u\in \cl{U}(k+j+1)\ |\ s_1(u)=u_i\}\ .\]
For the set of occupants of the frame $\theta$ we form a pullback
over $2^j$ sets as follows:

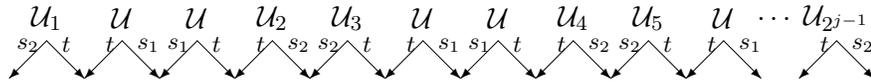
\begin{center}
\setlength{\unitlength}{1mm}
\begin{picture}(120,25)(0,-10)

\multiput(5,5)(10,0){10}{
\put(0,0){\vector(-1,-1){5}} %
\put(0,0){\vector(1,-1){5}}} %
\put(110,5){
\put(0,0){\vector(-1,-1){5}} %
\put(0,0){\vector(1,-1){5}}} %

\put(5,8){\makebox(0,0)[c]{$\cl{U}_1$}} %
\put(15,8){\makebox(0,0)[c]{$\cl{U}$}} %
\put(25,8){\makebox(0,0)[c]{$\cl{U}$}} %
\put(35,8){\makebox(0,0)[c]{$\cl{U}_2$}} %
\put(45,8){\makebox(0,0)[c]{$\cl{U}_3$}} %
\put(55,8){\makebox(0,0)[c]{$\cl{U}$}} %
\put(65,8){\makebox(0,0)[c]{$\cl{U}$}} %
\put(75,8){\makebox(0,0)[c]{$\cl{U}_4$}} %
\put(85,8){\makebox(0,0)[c]{$\cl{U}_5$}} %
\put(95,8){\makebox(0,0)[c]{$\cl{U}$}} %
\put(102,8){\makebox(0,0)[c]{$\dots$}}
\put(110,8){\makebox(0,0)[c]{$\cl{U}_{2^{j-1}}$}} %

\put(4,3.5){\makebox(0,0)[br]{{\footnotesize $s_2$}}}  %
\put(7,3.5){\makebox(0,0)[bl]{{\footnotesize $t$}}}

\put(14,3.5){\makebox(0,0)[br]{{\footnotesize $t$}}}  %
\put(17,3.5){\makebox(0,0)[bl]{{\footnotesize $s_1$}}}

\put(24,3.5){\makebox(0,0)[br]{{\footnotesize $s_1$}}}  %
\put(27,3.5){\makebox(0,0)[bl]{{\footnotesize $t$}}}

\put(34,3.5){\makebox(0,0)[br]{{\footnotesize $t$}}}  %
\put(37,3.5){\makebox(0,0)[bl]{{\footnotesize $s_2$}}}

\put(44,3.5){\makebox(0,0)[br]{{\footnotesize $s_2$}}}  %
\put(47,3.5){\makebox(0,0)[bl]{{\footnotesize $t$}}}

\put(54,3.5){\makebox(0,0)[br]{{\footnotesize $t$}}}  %
\put(57,3.5){\makebox(0,0)[bl]{{\footnotesize $s_1$}}}

\put(64,3.5){\makebox(0,0)[br]{{\footnotesize $s_1$}}}  %
\put(67,3.5){\makebox(0,0)[bl]{{\footnotesize $t$}}}

\put(74,3.5){\makebox(0,0)[br]{{\footnotesize $t$}}}  %
\put(77,3.5){\makebox(0,0)[bl]{{\footnotesize $s_2$}}}

\put(84,3.5){\makebox(0,0)[br]{{\footnotesize $s_2$}}}  %
\put(87,3.5){\makebox(0,0)[bl]{{\footnotesize $t$}}}

\put(94,3.5){\makebox(0,0)[br]{{\footnotesize $t$}}}  %
\put(97,3.5){\makebox(0,0)[bl]{{\footnotesize $s_1$}}}

\put(109,3.5){\makebox(0,0)[br]{{\footnotesize $t$}}}  %
\put(112,3.5){\makebox(0,0)[bl]{{\footnotesize $s_2$}}}

\end{picture}.
\end{center}

This completes the definition of $C_f$.

\section{Some examples at low dimensions}
\label{univdraw}

In this section we give some examples of elements of the
composition span $C = C_f$ for a 1-cell $f$.

\begin{itemize}
\item $j=0$
\end{itemize}

$C(0)$ is the set of universal 2-cells of the form

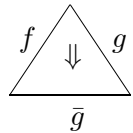
\sunit \begin{center}
\begin{picture}(18,30)(10,0)
\diagaw{10}{10}{1mm}
\end{picture} \end{center}

\noindent exhibiting $\bar{g}$ as a composite of $f$ and $g$.

\begin{itemize} \item $j=1$ \end{itemize}

We form a pullback over

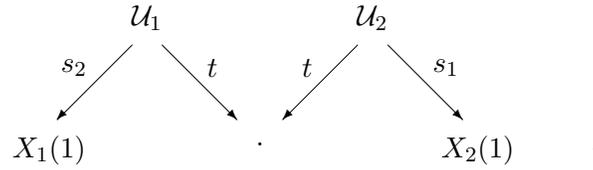
\sunit \begin{center}
\begin{picture}(65,40)(15,10)

\put(0,0){
\begin{picture}(80,30)     %
\put(18,23){\makebox(0,0)[t]{$X_1(1)$}}  
\put(31,37){\makebox(0,0)[b]{$\cl{U}_1$}}  
\put(46,22){\makebox(0,0)[t]{$.$}}  

\put(29,35){\vector(-1,-1){10}}  %
\put(33,35){\vector(1,-1){10}}  %
\put(23,32){\makebox(0,0)[r]{$s_2$}} 
\put(39,32){\makebox(0,0)[l]{$t$}}
\end{picture}}

\put(30,0){
\begin{picture}(80,30)     %
\put(18,23){\makebox(0,0)[t]{$$}}  
\put(31,37){\makebox(0,0)[b]{$\cl{U}_2$}}  
\put(45,23){\makebox(0,0)[t]{$X_2(1)$}}  

\put(29,35){\vector(-1,-1){10}}  %
\put(33,35){\vector(1,-1){10}}  %
\put(23,32){\makebox(0,0)[r]{$t$}} 
\put(39,32){\makebox(0,0)[l]{$s_1$}}
\put(60,21){.}
\end{picture}}
\end{picture}
\end{center}

\noindent So a typical element is of the form $(u_{31}, u_{32})$
with projections as shown below

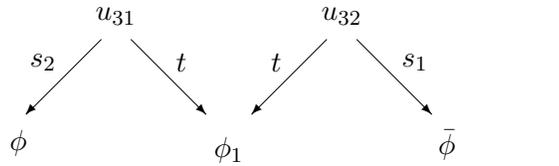
\begin{center}
\begin{picture}(65,40)(15,10)

\put(0,0){
\begin{picture}(80,30)     %
\put(18,23){\makebox(0,0)[t]{$\phi$}}  
\put(31,37){\makebox(0,0)[b]{$u_{31}$}}  
\put(46,22){\makebox(0,0)[t]{$\phi_1$}}  

\put(29,35){\vector(-1,-1){10}}  %
\put(33,35){\vector(1,-1){10}}  %
\put(23,32){\makebox(0,0)[r]{$s_2$}} 
\put(39,32){\makebox(0,0)[l]{$t$}}
\end{picture}}

\put(30,0){
\begin{picture}(80,30)     %
\put(18,23){\makebox(0,0)[t]{$$}}  
\put(31,37){\makebox(0,0)[b]{$u_{32}$}}  
\put(45,23){\makebox(0,0)[t]{$\bar{\phi}$}}  

\put(29,35){\vector(-1,-1){10}}  %
\put(33,35){\vector(1,-1){10}}  %
\put(23,32){\makebox(0,0)[r]{$t$}} 
\put(39,32){\makebox(0,0)[l]{$s_1$}}
\put(60,21){.}
\end{picture}}
\end{picture}
\end{center}

For example, the following two universal 3-cells exhibit
$\bar\phi$ as a composite of $f$ with $\phi$; this element of
$C(2)$ is in the frame $u_2 \lra u_2'$.

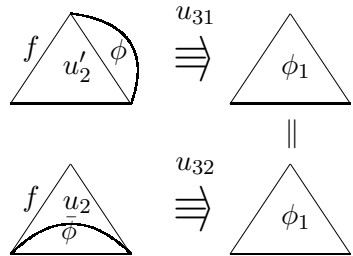
\begin{center}
\begin{picture}(45,60)(20,0) \diagaz{20}{15}{1mm} \end{picture}
\end{center}

\begin{itemize} \item $j=2$ \end{itemize}

We form a pullback over

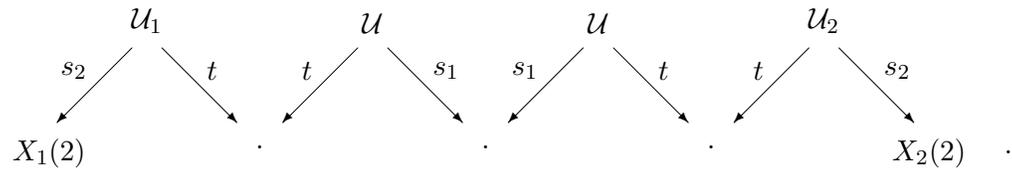
\sunit
\begin{picture}(80,40)(20,10)

\put(0,0){
\begin{picture}(80,30)     %
\put(18,23){\makebox(0,0)[t]{$X_1(2)$}}  
\put(31,37){\makebox(0,0)[b]{$\cl{U}_1$}}  
\put(46,22){\makebox(0,0)[t]{$.$}}  

\put(29,35){\vector(-1,-1){10}}  %
\put(33,35){\vector(1,-1){10}}  %
\put(23,32){\makebox(0,0)[r]{$s_2$}} 
\put(39,32){\makebox(0,0)[l]{$t$}}
\end{picture}}

\put(30,0){
\begin{picture}(80,30)     %
\put(18,23){\makebox(0,0)[t]{$$}}  
\put(31,37){\makebox(0,0)[b]{$\cl{U}$}}  
\put(46,22){\makebox(0,0)[t]{$.$}}  

\put(29,35){\vector(-1,-1){10}}  %
\put(33,35){\vector(1,-1){10}}  %
\put(23,32){\makebox(0,0)[r]{$t$}} 
\put(39,32){\makebox(0,0)[l]{$s_1$}}
\end{picture}}

\put(60,0){
\begin{picture}(80,30)     %
\put(18,23){\makebox(0,0)[t]{$$}}  
\put(31,37){\makebox(0,0)[b]{$\cl{U}$}}  
\put(46,22){\makebox(0,0)[t]{$.$}}  

\put(29,35){\vector(-1,-1){10}}  %
\put(33,35){\vector(1,-1){10}}  %
\put(23,32){\makebox(0,0)[r]{$s_1$}} 
\put(39,32){\makebox(0,0)[l]{$t$}}
\end{picture}}

\put(90,0){
\begin{picture}(80,30)     %
\put(18,23){\makebox(0,0)[t]{$$}}  
\put(31,37){\makebox(0,0)[b]{$\cl{U}_2$}}  
\put(45,23){\makebox(0,0)[t]{$X_2(2)$}}  

\put(29,35){\vector(-1,-1){10}}  %
\put(33,35){\vector(1,-1){10}}  %
\put(23,32){\makebox(0,0)[r]{$t$}} 
\put(39,32){\makebox(0,0)[l]{$s_2$}}
\put(55,21){.}
\end{picture}}

\end{picture}

A typical element is of the form $(u_{41},u_{42},u_{43},u_{44})$
with projections as shown below

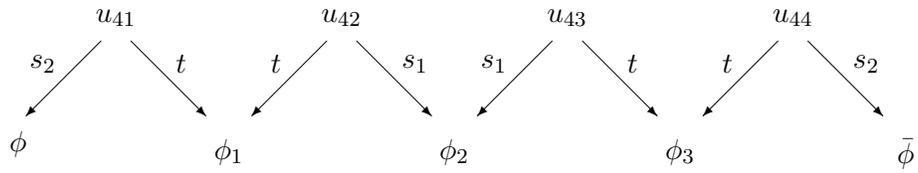
\begin{picture}(80,40)(20,10)

\put(0,0){
\begin{picture}(80,30)     %
\put(18,23){\makebox(0,0)[t]{$\phi$}}  
\put(31,37){\makebox(0,0)[b]{$u_{41}$}}  
\put(46,22){\makebox(0,0)[t]{$\phi_1$}}  

\put(29,35){\vector(-1,-1){10}}  %
\put(33,35){\vector(1,-1){10}}  %
\put(23,32){\makebox(0,0)[r]{$s_2$}} 
\put(39,32){\makebox(0,0)[l]{$t$}}
\end{picture}}

\put(30,0){
\begin{picture}(80,30)     %
\put(18,23){\makebox(0,0)[t]{$$}}  
\put(31,37){\makebox(0,0)[b]{$u_{42}$}}  
\put(46,22){\makebox(0,0)[t]{$\phi_2$}}  

\put(29,35){\vector(-1,-1){10}}  %
\put(33,35){\vector(1,-1){10}}  %
\put(23,32){\makebox(0,0)[r]{$t$}} 
\put(39,32){\makebox(0,0)[l]{$s_1$}}
\end{picture}}

\put(60,0){
\begin{picture}(80,30)     %
\put(18,23){\makebox(0,0)[t]{$$}}  
\put(31,37){\makebox(0,0)[b]{$u_{43}$}}  
\put(46,22){\makebox(0,0)[t]{$\phi_3$}}  

\put(29,35){\vector(-1,-1){10}}  %
\put(33,35){\vector(1,-1){10}}  %
\put(23,32){\makebox(0,0)[r]{$s_1$}} 
\put(39,32){\makebox(0,0)[l]{$t$}}
\end{picture}}

\put(90,0){
\begin{picture}(80,30)     %
\put(18,23){\makebox(0,0)[t]{$$}}  
\put(31,37){\makebox(0,0)[b]{$u_{44}$}}  
\put(46,22){\makebox(0,0)[t]{$\bar\phi$}}  

\put(29,35){\vector(-1,-1){10}}  %
\put(33,35){\vector(1,-1){10}}  %
\put(23,32){\makebox(0,0)[r]{$t$}} 
\put(39,32){\makebox(0,0)[l]{$s_2$}}
\end{picture}}

\end{picture}

\noindent exhibiting $\bar\phi$ as a composite of $f$ with $\phi$.
For example, the following element of $C(2)$ is in a frame with
target $(u_{31},u_{32})$.

\pagebreak

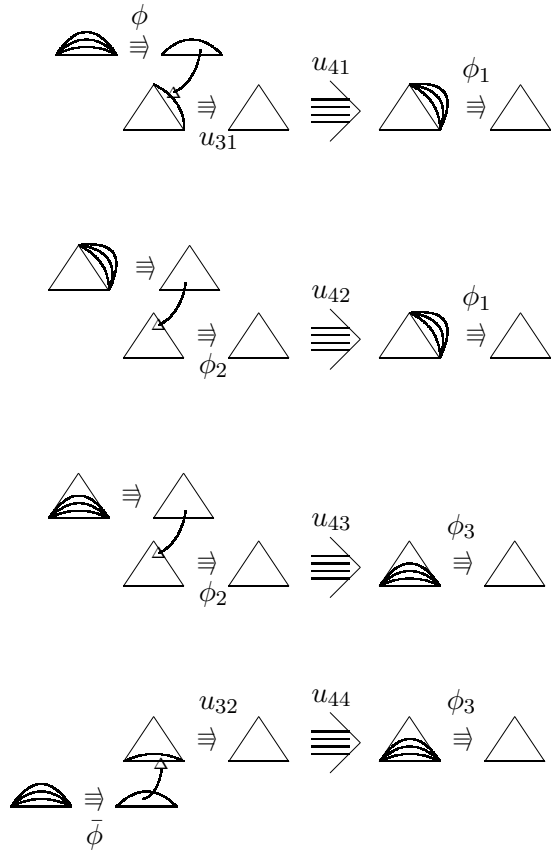
\sunit
\begin{center}
\begin{picture}(70,30)(10,-5) \diagba{20}{0}{1mm} \end{picture}

\begin{picture}(70,30)(10,-5) \diagbb{20}{0}{1mm} \end{picture}

\begin{picture}(70,30)(10,-5) \diagbc{20}{0}{1mm} \end{picture}

\begin{picture}(70,50)(10,-17) \diagbd{20}{15}{1mm} \end{picture}

\end{center}

\begin{itemize} \item $j=3$ \end{itemize}

Similarly, in $C(3)$ we have a typical element $(u_{51}, \ldots,
u_{58})$

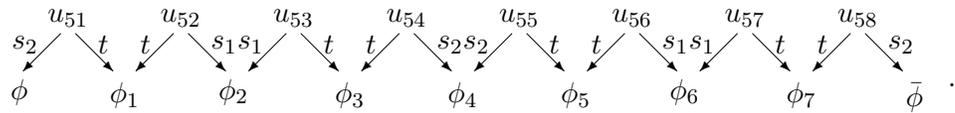
\setlength{\unitlength}{0.5mm}

\begin{picture}(80,40)(30,10)

\put(0,0){
\begin{picture}(80,30)     %
\put(18,23){\makebox(0,0)[t]{$\phi$}}  
\put(31,37){\makebox(0,0)[b]{$u_{51}$}}  
\put(46,22){\makebox(0,0)[t]{$\phi_1$}}  

\put(29,35){\vector(-1,-1){10}}  %
\put(33,35){\vector(1,-1){10}}  %
\put(23,32){\makebox(0,0)[r]{$s_2$}} 
\put(39,32){\makebox(0,0)[l]{$t$}}
\end{picture}}

\put(30,0){
\begin{picture}(80,30)     %
\put(18,23){\makebox(0,0)[t]{$$}}  
\put(31,37){\makebox(0,0)[b]{$u_{52}$}}  
\put(45,23){\makebox(0,0)[t]{$\phi_2$}}  

\put(29,35){\vector(-1,-1){10}}  %
\put(33,35){\vector(1,-1){10}}  %
\put(23,32){\makebox(0,0)[r]{$t$}} 
\put(39,32){\makebox(0,0)[l]{$s_1$}}
\end{picture}}

\put(60,0){
\begin{picture}(80,30)     %
\put(18,23){\makebox(0,0)[t]{$$}}  
\put(31,37){\makebox(0,0)[b]{$u_{53}$}}  
\put(46,22){\makebox(0,0)[t]{$\phi_3$}}  

\put(29,35){\vector(-1,-1){10}}  %
\put(33,35){\vector(1,-1){10}}  %
\put(23,32){\makebox(0,0)[r]{$s_1$}} 
\put(39,32){\makebox(0,0)[l]{$t$}}
\end{picture}}

\put(90,0){
\begin{picture}(80,30)     %
\put(18,23){\makebox(0,0)[t]{$$}}  
\put(31,37){\makebox(0,0)[b]{$u_{54}$}}  
\put(46,22){\makebox(0,0)[t]{$\phi_4$}}  

\put(29,35){\vector(-1,-1){10}}  %
\put(33,35){\vector(1,-1){10}}  %
\put(23,32){\makebox(0,0)[r]{$t$}} 
\put(39,32){\makebox(0,0)[l]{$s_2$}}
\end{picture}}

\put(120,0){
\begin{picture}(80,30)     %
\put(18,23){\makebox(0,0)[t]{$$}}  
\put(31,37){\makebox(0,0)[b]{$u_{55}$}}  
\put(46,22){\makebox(0,0)[t]{$\phi_5$}}  

\put(29,35){\vector(-1,-1){10}}  %
\put(33,35){\vector(1,-1){10}}  %
\put(23,32){\makebox(0,0)[r]{$s_2$}} 
\put(39,32){\makebox(0,0)[l]{$t$}}
\end{picture}}

\put(150,0){
\begin{picture}(80,30)     %
\put(18,23){\makebox(0,0)[t]{$$}}  
\put(31,37){\makebox(0,0)[b]{$u_{56}$}}  
\put(45,23){\makebox(0,0)[t]{$\phi_6$}}  

\put(29,35){\vector(-1,-1){10}}  %
\put(33,35){\vector(1,-1){10}}  %
\put(23,32){\makebox(0,0)[r]{$t$}} 
\put(39,32){\makebox(0,0)[l]{$s_1$}}
\end{picture}}

\put(180,0){
\begin{picture}(80,30)     %
\put(18,23){\makebox(0,0)[t]{$$}}  
\put(31,37){\makebox(0,0)[b]{$u_{57}$}}  
\put(46,22){\makebox(0,0)[t]{$\phi_7$}}  

\put(29,35){\vector(-1,-1){10}}  %
\put(33,35){\vector(1,-1){10}}  %
\put(23,32){\makebox(0,0)[r]{$s_1$}} 
\put(39,32){\makebox(0,0)[l]{$t$}}
\end{picture}}

\put(210,0){
\begin{picture}(80,30)     %
\put(18,23){\makebox(0,0)[t]{$$}}  
\put(31,37){\makebox(0,0)[b]{$u_{58}$}}  
\put(46,22){\makebox(0,0)[t]{$\bar\phi$}}  

\put(29,35){\vector(-1,-1){10}}  %
\put(33,35){\vector(1,-1){10}}  %
\put(23,32){\makebox(0,0)[r]{$t$}} 
\put(39,32){\makebox(0,0)[l]{$s_2$}}
\put(55,21){.} %
\end{picture}}

\end{picture}

\renewcommand{\qbeziermax}{100}

\noindent For example the following element of $C(3)$ (running
over two pages) has target $(u_{41},u_{42},u_{43},u_{44})$:

\newpage

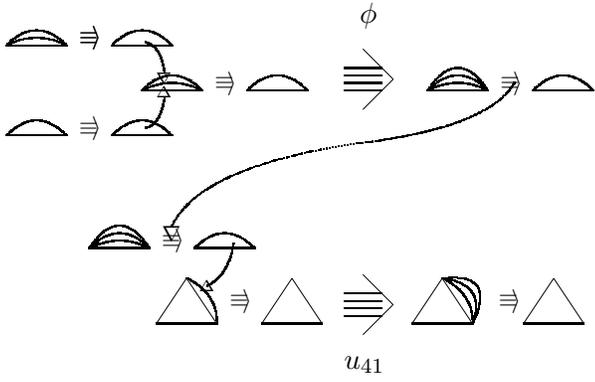
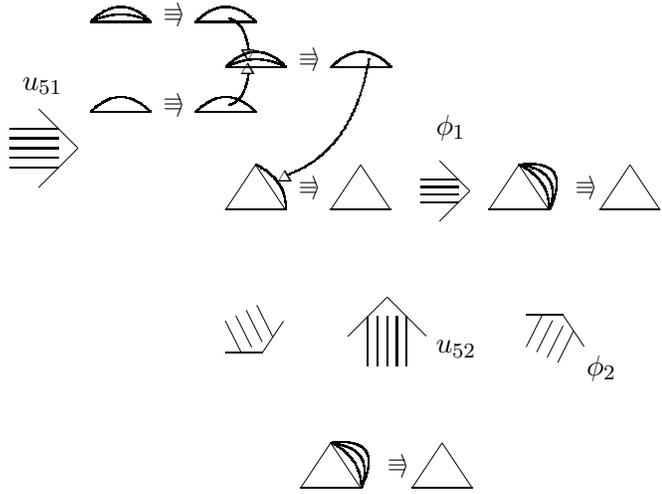
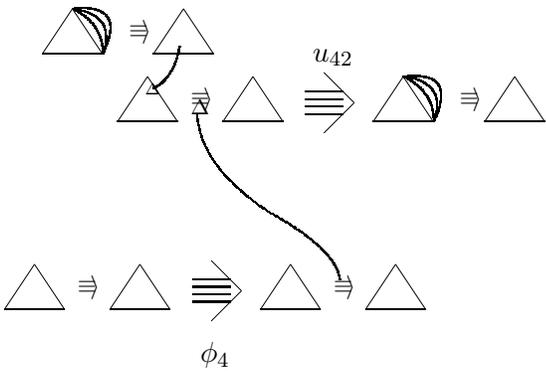
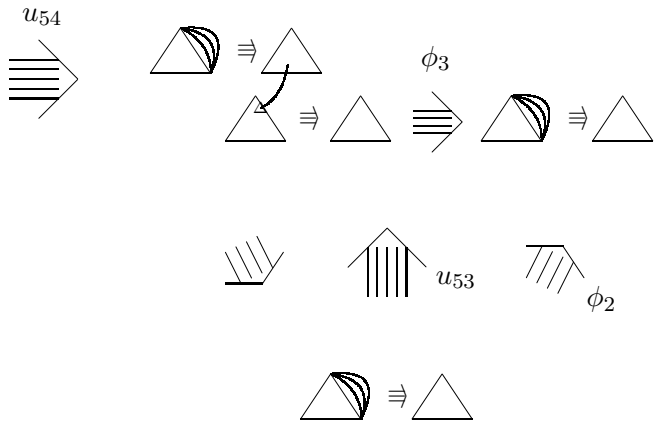
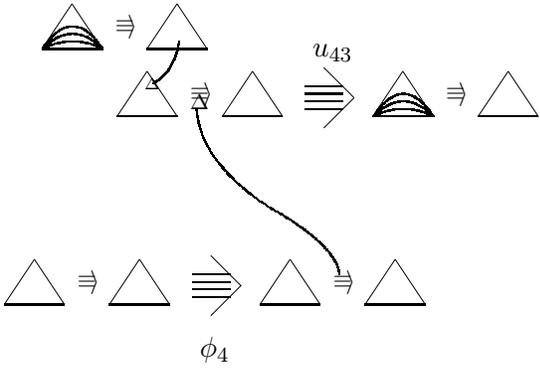
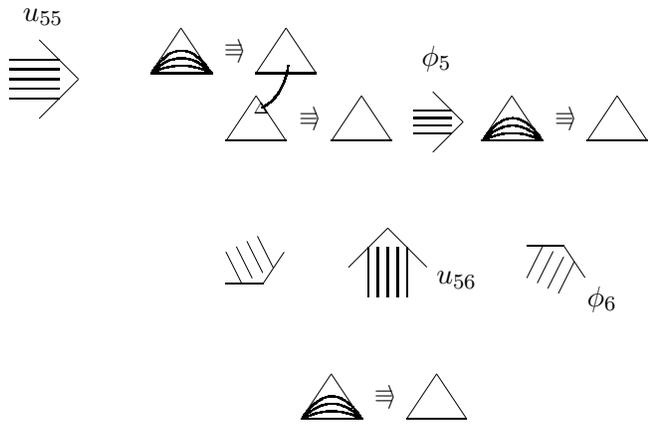
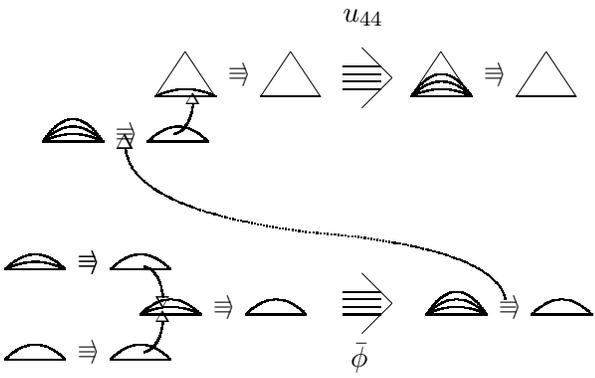
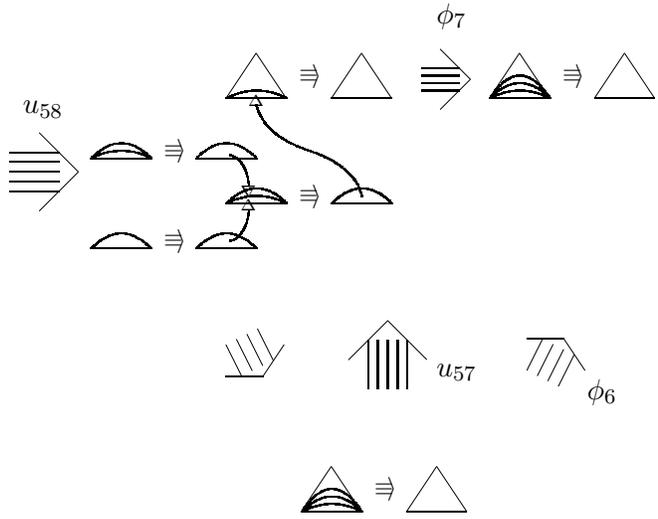
\sunit
\begin{picture}(20,80)(20,20) \diagbf{25}{35}{1mm} \end{picture}

\begin{picture}(20,80)(20,-20) \diagbg{25}{-15}{1mm} \end{picture}

\newpage

\begin{picture}(20,100) \diagbh{10}{35}{1mm} \end{picture}

\begin{picture}(20,100)(0,20) \diagbi{10}{60}{1mm} \end{picture}

\section{Comparison}
\label{comp}

We now compare the new characterisation with the original definition for $n
\leq 2$, and show that at these low dimensions the notions do indeed coincide. 
We argue explicitly; at higher dimensions such an approach rapidly becomes
unfeasible.  An effective method for handling such algebra is urgently needed.  

For convenience we refer to the `old' and `new' universal properties as Property
1 (P1) and Property 2 (P2) respectively, and we continue to use all earlier notation.

\bpoint{$n=0$}

P1 and P2 clearly coincide for $k$-cells with $k >n$; this deals with all
possibilities when $n=0$.

\bpoint{$n=1$}

Let $X$ be an opetopic 1-category.  We show that P1 and P2 are equivalent for
$k=1$.  Let $f$ be a 1-cell $a \map{f} b$ in $X$.  $f$ has P1 if and only if for
all $a \map{g} c \in X$ there is a unique factorisation

\begin{picture}(25,22)(-3.5,-4)
\twotwo{$$}{$$}{$$}{$f$}{$\bar{g}$}{$g$}{$u$}
\end{picture}

\noindent (see \cite{che10}) i.e. $\tau_0$ is an isomorphism $C(0) \lra X_2(0)$,
which says precisely that $f$ has P2.  

Note that this argument immediately generalises to all $n=k$.  

\bpoint{$n=2$}

Let $X$ be an opetopic 2-category.  The above arguments deal with $k \geq 2$;
we show that the notions coincide when $k=1$.  

Let $a \map{f} b$ have P1.  We show that 
	\[\tau_f: C_f \lra X(a,c)\]
is an equivalence of 1-categories.

\begin{enumerate}

\item We show that it is essentially surjective on objects; in fact it is
surjective `on the nose'.  A 0-cell in the codomain is a 1-cell of the form $a
\map{g} c$.  Given any such, there certainly is a factorisation

\begin{center}
\begin{picture}(25,22)(-3.5,-4)
\twotwo{$$}{$$}{$$}{$f$}{$\bar{g}$}{$g$}{$u$}
\end{picture}
\end{center}

so we have $\bar{g}$ such that
	\[\tau_0: \bar{g} \mapsto g.\]

\item We show that is is `locally an isomorphism'.  Consider the 1-frame in
$C_f$

\begin{center}
\begin{picture}(25,12)(-3.5,2)
\twotwo{$$}{$$}{$$}{$f$}{$\bar{g}_1$}{$g_1$}{$u_1$}
\end{picture}
\hspace{1cm} \map{?} \hspace{1cm}
\begin{picture}(25,12)(-3.5,2)
\twotwo{$$}{$$}{$$}{$f$}{$\bar{g}_2$}{$g_2$}{$u_2$}
\end{picture}.
\end{center}
\vspace{3ex}
We show that $\tau$ is an isomorphism on the 0-category (set) induced by this
frame.  So consider 1-cell in the codomain
	
\begin{center}
\begin{picture}(35,15)(5,-3)
\onetwo{10}{0}{1mm}{$$}{$$}{$g_1$}{$g_2$}{$\alpha$}
\end{picture}.
\end{center}
We have

\begin{center}
\begin{picture}(25,22)(-5,-4)
\diagbr{}{}{}{}{}{}{}{$u_1$}{$\alpha$}
\end{picture}
\begin{picture}(10,20)
\three{5}{10}{1mm}{$$}
\end{picture}
\begin{picture}(25,22)(-5,-4)
\twotwob{$$}{$$}{$$}{$$}{$$}{$$}{$\theta$}
\end{picture},
\end{center}
say, and by P1, $\theta$ induces a unique factorisation
\begin{center}
\begin{picture}(25,22)(-5,-4)
\diagbt{}{}{}{}{}{}{}{$\bar{\alpha}$}{$u_2$}
\end{picture}
\begin{picture}(10,20)
\three{5}{10}{1mm}{$$}
\end{picture}
\begin{picture}(25,22)(-5,-4)
\twotwob{$$}{$$}{$$}{$$}{$$}{$$}{$\theta$}
\end{picture}
\end{center}
so we have a unique pre-image of $\alpha$ as required.

\end{enumerate}

The converse follows easily.

\section{Conclusions}
\label{conc}

We conclude that although the outline of the basic syntax of opetopes seems
secure, universality is less well understood, and we remain unsure of the
ideal form in which it should be defined. The alternative characterisation
described in this work seems right in `spirit', but in the end
the mathematics that emerges is not as `slick' as might be hoped.
It therefore appears that there is much scope for further work in this
area.

\addcontentsline{toc}{section}{References}
\bibliography{bib0209}

\end{document}